\documentclass[10pt]{article}
\usepackage{inputenc}
\usepackage[english]{babel}
\usepackage{amssymb,amsmath,amsthm,amscd}
\usepackage{color}   

\newfont{\rams}{msbm10 scaled\magstep1}
\newfont{\ramss}{msbm10 scaled\magstep0}
\newfont{\iams}{msbm10}
\newfont{\gotic}{eufm10 scaled\magstep1}
\newfont{\bellap}{eusm10 scaled\magstep1}

\newcommand{\nat}{\mbox{\rams \symbol{'116}}}

\newcommand{\q}{{J}}
\newcommand{\p}{{P}}

\newcommand{\da}{{A\!\Join^f\!\!J}}

\newcommand{\sss}{{\rm Spec}}
\newcommand{\Jac}{{\rm Jac}}%
\newcommand{\hgt}{{\rm ht}}
\newcommand{\mmm}{{\rm Max}}

\newcommand{\Ker}{{\rm Ker}}
\newcommand{\tot}{{\rm Tot}}

\newcommand{\z}{{\ldots}}
\newcommand{\w}{{\setminus}}

\newcommand{\ude}{\mbox{\textsl{id}}} %

\swapnumbers
\newtheoremstyle{break}
  {9pt}
  {9pt}
  {\itshape}
  {}
  {\textsc}
  {.}
  {.7em}
  {}
\newtheoremstyle{break1}
  {9pt}
  {9pt}
  {}
  {}
  {\textsc}
  {.}
  {.7em}
  {}
\theoremstyle{break}
\newtheorem{thm}{ \textsc{Theorem}}[section]
\newtheorem{defn}[thm]{ \textsc{Definition}}
\newtheorem{cor}[thm]{ \textsc{Corollary}}
\newtheorem{lem}[thm]{ \textsc{Lemma}}
\newtheorem{prop}[thm]{ \textsc{Proposition}}
\theoremstyle{break1}
\newtheorem{ex}[thm]{ \textsc{Example}}
\newtheorem{oss}[thm]{ \textsc{Remark}}
\theoremstyle{remark}

\title{Properties of chains of prime ideals in an  amalgamated algebra along an ideal 
}
\author{Marco D'Anna \footnote{Partially supported by a Grant MIUR-PRIN.}
\and Carmelo A. Finocchiaro
\and Marco Fontana}

\begin{document}

\maketitle

\hfill{\footnotesize \sl }%

\bigskip

\begin{abstract}
Let $f:A \rightarrow B$ be a ring homomorphism and let $J$ be an ideal of $B$. In this paper, we study the
amalgamation of $A$ with $B$ along $J$ with respect to $f$ (denoted by $\da$), a construction that provides a
general frame for studying  the amalgamated duplication of a ring   along an ideal, introduced and studied by
D'Anna and Fontana in 2007, and  other classical constructions (such as the $A+ XB[X]$, the $A+ XB[\![X]\!]$ and
the $D+M$ constructions). In particular, we completely describe the prime spectrum of the amalgamated duplication
and we give bounds for its Krull dimension.
\end{abstract}

\bigskip

\section{Introduction}

Let $A$ and $B$ be commutative rings with unity, let $J$ be  an ideal of $B$ and let $f:A\longrightarrow B$ be a
ring homomorphism. In this setting, we can consider the following subring of $A\times B$:
$$
A\Join^f\! J := \{(a,f(a)+j) \mid  a \in A, \ j \in J \}
$$called {\it the amalgamation of $A$ with $B$ along $J$ with respect to $f$}. This construction is a generalization
of the amalgamated duplication of a ring along an ideal (introduced and studied in \cite{d'a-f-1}, \cite{d'a},
\cite{d'a-f-2} and in \cite{m-y}).  Moreover, several classical constructions (such as the $A+XB[X]$, the
$A+XB[\![X]\!]$ and the $D+M$ constructions) can be studied as particular cases of the amalgamation \cite[Examples
2.5 and 2.6]{dafifoproc} and other classical constructions, such as the Nagata's idealization (cf. \cite[page
2]{N}, \cite[Chapter VI, Section 25]{H}), also called Fossum's trivial extension (cf. \cite{Fos}, \cite{K-M} and
\cite{B-K-M}), and the CPI extensions (in the sense of Boisen and Sheldon \cite{bo}) are strictly related to it
\cite[Example 2.7 and Remark 2.8]{dafifoproc}.

On the other hand, the amalgamation $\da$ is related to a construction proposed by D.D. Anderson in \cite{a-06}
and motivated by a classical construction due to Dorroh \cite{do1}, concerning the embedding of a ring without
identity in a ring with identity. An ample introduction on the genesis  of the notion of amalgamation is given in
\cite[Section 2]{dafifoproc}.

One of the key tools for studying $\da$ is based on the fact that the amalgamation  can be studied in the frame of
pullback constructions   \cite[Section 4]{dafifoproc}  (for a systematic study of this type of  constructions, cf.
\cite{F}). This point of view allows us to deepen the study initiated in \cite{dafifoproc} and to provide an ample
description of various properties of $\da$, in connection with the properties of $A$, $J$ and $f$.

 More precisely,
in  \cite{dafifoproc},  we studied the basic properties of this construction (e.g., we provided characterizations
for $\da$ to be a Noetherian ring, an integral domain, a reduced ring)
  and we characterized those distinguished pullbacks that can be expressed as an
amalgamation.
In this paper, we pursue
the investigation on the structure of the rings of the form $\da$, with particular attention to the  prime
spectrum, to the chain properties and to the Krull dimension. In particular, after recalling (in Section 2) some
basic properties proved in \cite{dafifoproc} and \cite{dafifoaltro}, needed in the present paper,
we start our
investigation by deepening the study of chains of prime ideals in  pullback constructions (Proposition
\ref{dimp}).

 In Section 3, we study the integral closure of $\da$ in its total ring of fractions and, finally, in Section 4,
we concentrate our attention to evaluate its Krull dimension. In particular, we provide upper and lower
bounds for dim($\da$) (Proposition \ref{lower} and Theorem \ref{carinissima}) and we show that these bounds, obtained in a
such  general setting, are so sharp that
 generalize, and possibly improve,
analogous bounds established for the very particular cases of integral domains of the form $A+XB[X]$ \cite{fik} or
$A+XB[\![X]\!]$ \cite{doka}.


\section{Preliminaries}
Before beginning a systematic study of the ring $\da$, we recall from our introductory paper \cite{dafifoproc} to
the subject some basic properties of this construction.
\begin{prop}
{\rm \cite[Proposition 5.1]{dafifoproc}}\label{inizio}   Let $f: A\rightarrow B$ be a ring  homomorphism,  $\q$ an ideal
of $B$  and set $\da := \{ (a, f(a)+j) \mid a\in A, \ j \in J \}$.
\begin{enumerate}
 \item[\rm (1)]  Let  $\iota := \iota_{A, f, J}: A\rightarrow \da$ be the natural ring homomorphism defined by $\iota(a)  := (a, f(a))$, for all $a \in A$. Then $\iota$ is an embedding, making
$\da$   a  ring   extension of $A$ \  (with $\iota(A) =  \Gamma(f) \  (:=\{(a,f(a)) \mid a\in A\}$  subring of
$\da$).

\item[\rm (2)]  Let $I$ be an ideal of $A$ and set $ I\!\Join^f\!\! J :=\{(i, f(i)+j) \mid i\in I, j \in J \}$.
Then $I\!\Join^f\!\! J$  is an ideal of $\da$, the composition  of canonical homomorphisms
$A\stackrel{\iota}{\hookrightarrow} \da\twoheadrightarrow \da/I\!\Join^f\!\! J$ is a surjective ring homomorphism
and its kernel coincides with  $I$.\\ Hence, we have the following canonical isomorphism:
 $$
 \frac{A}{I}\cong\frac{\da}{I\!\Join^f\!\! J}\,.
$$

 \item[\rm (3)] Let $p_{_A}: \da \rightarrow A$ and $p_{_B}:\da\rightarrow B$ be the natural projections of
 $\da \subseteq A \times B$ {into} $A$ and $B$, respectively.  Then $p_{_A}$ is surjective
 and\, $\Ker(p_{_A})=\{0\}\times J$.\\
Moreover, $p_{_B}(\da)=f(A)+ J$ and $\Ker(p_{_B})=f^{-1}(\q)\times \{0\}$.  Hence, the following canonical
isomorphisms  hold:
$$
\frac{\da}{(\{0\}\times \q)}\cong A\quad \mbox{ and } \quad\frac{\da}{f^{-1}(\q)\times \{0\}}\cong f(A)+\q \,.
$$
\item[\rm(4)] Let $\gamma:\da\rightarrow (f(A)+J)/J$
be the natural ring homomorphism, defined by $(a,f(a)+j)\mapsto f(a)+J$. Then $\gamma$ is surjective and
$\Ker(\gamma)=f^{-1}(J)\times J$. Thus, there exists a natural isomorphism
$$
\frac{\da}{f^{-1}(J)\times J}\cong\frac{f(A)+J}{J}\,.
$$
In particular, when $f$ is surjective, we have the following natural isomorphism
$$
\frac{\da}{f^{-1}(J)\times J}\cong \frac{B}{J}\,.
$$
 \end{enumerate}
\end{prop}
 Recall
that, in \cite{dafifoproc}  and \cite{dafifoaltro}, we have shown that the ring $\da$ can be re\-presented as a
pullback of natural ring homomorphisms and, using the notion of ring retraction,  we have characterized the
pullbacks that produce exactly rings of the form $\da$ (see also Propositions \ref{pull} and
\ref{ringretraction}).
 Now we will make some pertinent
remarks and prove a new result on chains of prime ideals of pullbacks, that will be useful for our subsequent
investigation of the ring $\da$.
\medskip
\begin{defn}\label{de} \rm
We recall that, if $\alpha:A\rightarrow C,\,\,\, \beta:B\rightarrow C$ are ring homomorphisms, the subring
$D:=\alpha\! \times_{_C}\! \beta:=\{(a,b)\in A\times B \mid \alpha(a)=\beta(b)\}$ of $A\times B$ is called the
\textit{{pullback}} (or \textit{fiber product}) of $\alpha$ and $\beta$.
%
 In the following, we will denote by   $p_{_A}$ (respectively, $p_{_B}$)  the restriction to $\alpha\! \times_{_C}\! \beta$ of the projection of
$A\times B$  onto $A$ (respectively, $B$).
\end{defn}

The following proposition is a straightforward consequence of the definitions.

\begin{prop}{\rm \cite[ Proposition 4.2]{dafifoproc}\label{pull}}
Let $f:A\rightarrow B$ be a ring homomorphism and $\q$ be an ideal of $B$. If \ $\pi:B\rightarrow B/J$ is the
canonical projection and $\breve f:=\pi\circ f$, then $\da=\breve f\times_{_{B/J}}\pi$.
\end{prop}
\noindent Now, recall that a ring homomorphism $r:B\rightarrow A$ is called {\it a ring retraction} if there
exists an (injective) ring homomorphism $i:A\rightarrow B$ such that $r\circ i=\ude_A$. In this case, we say
also that $A$ is {\it a retract of} $B$.
\begin{ex}{\cite[Remark 4.6]{dafifoproc}}\label{retrex}
Let $f:A\rightarrow B$ be a ring homomorphism and $J$ an ideal of $B$.
Then $A$ is a retract of $\da$ and the map $p_A:\da\rightarrow A$, defined in Proposition \ref{inizio}(3), is a ring
retraction. In fact, we have $p_A\circ \iota= \ude_A$, where $\iota$ is the ring embedding of $A$ into $\da$
(Proposition \ref{inizio}(1)).
\end{ex}
The pullbacks of the form $\da$ form a distinguished subclass of
the class of pullbacks of ring homomorphisms, as described in the following proposition.
\begin{prop} {\rm \cite[Proposition 4.7]{dafifoproc}}\label{ringretraction}
Let $A,B,C,\alpha,\beta,p_{_A},p_{_B}$ be as in Definition \ref{de}. Then, the following conditions are
equivalent.
\begin{enumerate}
\item[\rm(i)] 
  $p_{_A}:\alpha\! \times_{_C}\! \beta\rightarrow A$ is a ring retraction.
\item[\rm(ii)] There exist an ideal $J$ of $B$ and a ring homomorphism
$f:A\rightarrow B$ such that $\alpha\! \times_{_C}\! \beta=\da$. \hfill $\Box$
\end{enumerate}\end{prop}
  Let  $f: A\rightarrow B$ be a ring homomorphism, and set $X:=\sss(A), \ Y:=\sss(B)$. Recall that
$f^*:Y\rightarrow X$ denotes the continuous map  (with respect to the Zariski topologies)  naturally associated to
$f$  (i.e.,  $f^*(Q):= f^{-1}(Q)$ for all $Q\in Y$).  Let  $S$ be a subset of $A$. Then, as usual,   $V_X(S)$, or
simply $V(S)$, if no confusion can arise, denotes the closed subspace of $X$, consisting of all prime ideals of
$A$ containing $S$.  We will denote by ${\rm Jac}(A)$ the Jacobson radical of a ring $A$ and we will call \it
local ring \rm a (not necessarily Noetherian) ring with a unique maximal ideal.


Now, we collect some results about the structure of the prime ideals of the ring $\da$. The proof of the following
proposition is based on well known properties of rings arising from pullbacks \cite[Theorem 1.4]{F} (for details,
see \cite {dafifoaltro}).

\begin{prop}\label{spec}
 With the notation of Proposition \ref{inizio},  set $X:=\sss(A)$, $Y:=\sss(B)$,  and $W:=\sss(\da)$,
 and $J_0 :=\{0\}\times J \ (\subseteq \da)$. For all $P \in X$ and $Q \in Y$, set:
$$
\begin{array}{rl}
P^{\prime_{_{\! f}}}:= &\hskip -4pt  P\!\Join^f\!\! J :=\{(p, f(p)+j ) \mid p\in \p,\ j\in J \}\,,\\
\overline Q^{_{_{ f}}}:=&\hskip -4pt \{(a, f(a)+j) \mid a\in A,\  j \in J,\ f(a)+ j \in Q\}\,.
\end{array}
$$
 Then,  the following statements hold.
\begin{enumerate}
\item [\rm (1)] The map $P \mapsto P^{\prime_{_{\! f}}}$ establishes a closed embedding of $X$ into  $W$,
 so its image, which coincides with $V(J_0)$, is homeomorphic to $X$.
\item [\rm (2)] The map $Q \mapsto\overline Q^{_{_{f}}} $ is {{a}} homeomorphism of $Y\w V(J)$ onto
$W\w V(J_0)$.

\item [\rm (3)] The prime ideals of $\da$ are of the type\ $P^{\prime_{_{\! f}}}$ or \ $\overline Q^{_{_{ f}}}$\!,\,  for
 $P$ varying in $ X$ and $Q$ in $Y\w V(\q)$.
\item [\rm(4)] Let $P \in\sss(A)$.  Then,  $P^{\prime_{_{\! f}}}$
 is a maximal ideal of $\da$ if and only if $P$ is a maximal ideal of $A$.
  \item [\rm(5)] Let $Q$ be a prime ideal of $B$ not containing $J$. Then,
   $\overline Q^{_{_{f}}}$ is a maximal ideal of $\da$ if and only if $Q$ is a maximal ideal of $B$.\end{enumerate}
In particular:
$$
\mmm(\da)=\{P^{{\prime_{_{\! f}}}} \mid  P \in\mmm(A)\}  \cup \{\overline Q^{_{_{f}}} \mid Q \in \mmm(B)\w V(\q)\}\,.
$$
\vskip -20pt\hfill $\Box$
\end{prop}
\medskip
The last result of this section concerns the chains of prime ideals in rings arising from pullbacks of rather general type.
\medskip

\begin{prop}\label{dimp}
With the notation of Definition \ref{de}, assume $\beta$  surjective. Let $H'$ and $H''$ be prime ideals of $D$
such that $H'\subsetneq H''$. Assume that  $H'\in\sss(D)\w$ $ V(\Ker(p_{_A}))$,  $H '' \in V(\Ker(p_{_A}))$, and
that $H'$ and $H''$ are adjacent prime ideals.  Then, there exist  two prime ideals $Q'$ and $Q''$ of $B$, with  $Q' \subsetneq Q''$, and moreover such that $Q' \notin
V(\Ker(\beta))$, $p_B^{-1}(Q') = H'$, and $p_{_B}^{-1}(Q'')=H''$.
\end{prop}
\noindent \textsc{Proof}.    Note that the existence (and uniqueness) of a prime ideal $Q'$ of $B$ such that $Q' \notin
V(\Ker(\beta))$ and $p_B^{-1}(Q') = H'$  is well known \cite[Theorem 1.4, Statement (c) of the proof]{F}.

On the other hand, note that $p_{_B}^{-1}(L+\Ker(\beta))=p_{_B}^{-1}(L)+\Ker(p_{_A})$, for each ideal $L$ of $B$. Now,
it is clear that the set
$$
{\boldsymbol{\mathcal{S}}}(Q'):=\{L \mbox{ ideal of }B \mid  Q'+\Ker(\beta)\subseteq L \mbox { and  }
p_{_B}^{-1}(L)\subseteq H''\}
$$
is nonempty (it contains $Q'+\Ker(\beta)$)  and inductive. Thus, by Zorn's lemma, ${\boldsymbol{\mathcal{S}}}(Q')$
contains a maximal element $Q''$, which is easy to see that is a prime ideal of $B$. Since $H'' \supseteq
p_{_B}^{-1}(Q'')\supseteq p_{_B}^{-1}(Q')+\Ker(p_{_A})\supsetneq H'$ and $H',H''$ are adjacent prime ideals, we
have $p_{_B}^{-1}(Q'')= H''$.\hfill$\Box$

 \section{Integral closure of the ring $\da$} Given  a ring extension $R\subseteq S$,   the integral closure
of $R$ in $S$ will be denoted by $\overline{R}^{\mbox{\tiny \it S}}$; the integral closure of $R$ in its total
ring of fractions ${\rm Tot}(R)$
  will be simply denoted by $\overline{R}$.

Now, we want to determine the integral closure of the ring $\da$ in its total ring of fractions. It is easy to
compute $\tot(\da)$ in some cases.
\begin{prop}\label{tot}  Let $f: A\rightarrow B$ be a ring
homomorphism,  $\q$ an ideal of $B$,  and let $\da$ be as in  Proposition \ref{inizio}.
Assume that  $\q$ and $f^{-1}(\q)$ are regular ideals of $B$ and $A$, respectively.
Then $\tot(\da)$ is canonically isomorphic to $ \tot(A)\times \tot(B)$.

\end{prop}
\noindent
\textsc{Proof.}
Note that $J_1:=f^{-1}(J)\times J$ is the conductor of $\da$ in  $A\times B$ (i.e., the largest  ideal of $\da$ that is also an ideal of $A\times  B$). Since both $f^{-1}(J)$ and $J$ are regular ideals, then  $J_1$ is a regular ideal of $A\times B$. Now, the conclusion follows immediately by applying \cite[pag. 326]{g}.\hfill$\Box$

\begin{oss}
Note that, in Proposition \ref{tot}, the assumption that $\q$ and $f^{-1}(\q)$ are regular ideals is essential.
For example,  let $A$ be an integral domain with quotient field $K$,  $B$  an overring of $A$,  and  let
$\q=\{0\}$.  Then, in this situation,  $  \da \cong A$ (Proposition \ref{pull}), and thus $\tot(\da)$ is
isomorphic to $K$, but $\tot(A)\times\tot(B)=K\times K$.

In the previous example, $J$ and $f^{-1}(J)$ are both the zero ideal. Another example, for which  $J$  is a
nonzero regular ideal,  is given next.  Let  $A$ be an integral domain with quotient field $K$, set $B:=A[X]$ and
$\q:=(X)$, and let $f: A \hookrightarrow A[X]$ be the natural inclusion.  In this case, from Proposition
\ref{pull} we deduce that  $\da \cong A + XA[X]=A[X]$,  and hence $\tot(\da)=K(X)$. However,
$\tot(A)\times\tot(B)=K\times K(X)$. (Note that in this example $f^{-1}(\q)= A \cap \q=\{0\}$.)

Another example, for which both $J$ and $f^{-1}(J)$ are nonzero and not regular ideals, is the following. Let $K$
be a field  and set $A:=K^{(3)}, B:=K^{(2)}$,  and $J:=\{0\}\times K$, where $K^{(n)}$ is the direct product  ring
$K\times K\times ... \times K$ ($n$--times).   If $f$ is the projection defined by $(a,b,c)\mapsto (a,b)$, it is
immediately seen that $\da\cong K^{(4)}$. Then $\tot(\da)\cong K^{(4)}$, but $\tot(A)\times \tot(B)\cong K^{(5)}$.
\end{oss}
\noindent

We have already observed in \cite[Section 5]{dafifoproc} that the ring $B_{\diamond} := f(A) +J$
(subring of $B$) plays a relevant role in the
construction $\da$.  The next result provides further evidence to this fact.

\begin{lem}\label{integrale} Let $f: A\rightarrow B$ be a ring homomorphism,
$\q$ an ideal of $B$, and let $\da$ be as in  Proposition \ref{inizio}.
The ring $A\times (f(A)+\q)$, subring of $A \times B$, which contains $\da$   is integral over $\da$. More
precisely, every element of $A\times (f(A)+\q)$ has degree  at most two   over $\da$.
\end{lem}
\noindent \textsc{Proof}.  Let $(\alpha, f(a)+j)\in A\times (f(A)+\q)$ with $\alpha, a \in A$ and $j \in J$.
Assume that $(\alpha, f(a)+j) \notin \da$, thus, in particular, $\alpha \neq a$.  Then,  the element $(\alpha,
f(a)+j)$ is a root of the monic polynomial $(X-(\alpha, f(\alpha)))(X-(a, f(a)+j))\in (\da)[X]$. \hfill$\Box$
\begin{prop} \label{chius}  With the notation of Lemma \ref{integrale}, assume that $J$ and $f^{-1}(\q)$
are regular ideals of $B$ and $A$, respectively.
Then  $\overline{\da}$  (i.e., the integral closure of $\da$ in its total ring of fractions) coincides with $\overline{A}
\times \overline{f(A)+\q}$. In particular, if $f$ is an integral homomorphism, then $\overline{\da}
=\overline{A} \times\overline{B}$.
\end{prop}
\noindent \textsc{Proof}.  Recall that, under the present hypothesis on $J$ and $f^{-1}(\q)$, we have ${\rm
Tot}(\da)={\rm Tot}(A \times B)$,  which is canonically isomorphic to ${\rm Tot}(A) \times {\rm Tot}(B)$
(Proposition \ref{tot}).  Therefore, it is easy to see that $\overline{\da} \subseteq \overline{A} \times
\overline{f(A)+\q}$. On the other hand, the ring $\overline{A} \times \overline{f(A)+\q}$ is obviously integral
over $A \times (f(A)+\q)$ and $A \times (f(A)+\q)$ is integral over $\da$ (Lemma \ref{integrale}). Thus
$\overline{A} \times \overline{f(A)+\q}$ is integral over $\da$. The conclusion is now straightforward.
\hfill$\Box$
\bigskip

\begin{oss} If  we do not assume that $J$ and $f^{-1}(\q)$ are regular ideals of $B$
and $A$, respectively, then the argument used in the proof of Proposition \ref{chius} shows that the integral
closure of $\da$  in  ${\tot(A)\times \tot(B)}$ coincides with $\overline{A} \times \overline{f(A)+\q}$.
\end{oss}
Now, we want to investigate when the ring $\da$ is integral over $\Gamma(f)(:=\{(a,f(a))\mid a\in A\})$.
\begin{lem}
\label{integralegamma} Let $f:A\longrightarrow B$, $\q\subseteq B$, and $\da$ be as in  Proposition \ref{inizio}.
Then, the following conditions are equivalent.
\begin{enumerate}
\item [\rm (i)] $f(A)+\q$ is integral over $f(A)$.
\item [\rm (ii)] $\da$ is integral over $\Gamma(f)$.
\end{enumerate}
In  particular, if $f$ is an integral homomorphism, then $\da$ is integral over $\Gamma(f)  \ (\cong A)$.
\end{lem}
\noindent \textsc{Proof}. (i) implies (ii). Let $(a,f(a)+j)$ be a  nonzero element of $\da$. Thus, by condition
(i), there exist a positive integer $n$ and $a_0, a_1,\z,a_{n-1}\in A$ such that
$(f(a)+j)^n+\sum_{i=0}^{n-1}f(a_i)(f(a)+j)^i=0$. Therefore, it is easy to verify that $(a,f(a)+j)$ is a root of
the monic polynomial $[X-(a,f(a))][X^n+\sum_{i=0}^{n-1}(a_i,f(a_i))X^i]\in \Gamma(f)[X]$. Conversely, consider an
element $f(a)+j\in f(A)+\q$. By condition (ii), $(a,f(a)+j)$ is integral over $\Gamma(f)$, and hence the equation
of integral dependence of $(a,f(a)+j)$ over $\Gamma(f)$ gives us the equation of integral dependence of $f(a)+j$
over $f(A)$. The last statement is straightforward.\hfill $\Box$

\section{Krull dimension of $\da$}

Now, we want to  study   the   Krull   dimension of the ring $\da$.   We start with an easy observation.

\begin{prop}\label{su}
 Let $f: A\rightarrow B$,  $\q$,  and   $\da$ be as in  Proposition \ref{inizio}. Then
$\dim(\da)=\max\{\dim(A),\dim(f(A)+\q)\}$.  In particular, if $f$ is surjective,
then $\dim(\da)=\max\{\dim(A),\dim(B)\}=\dim(A)$.
\end{prop}
\noindent \textsc{Proof}. By Lemma \ref{integrale} and \cite[Theorem 48]{k}, it follows immediately that
$\dim(\da)=\dim(A\times(f(A)+\q))$. Thus, the  conclusion is an easy  consequence of the fact that
$\sss(A\times(f(A)+\q))$ is  canonically  homeomorphic to the disjoint union of $\sss(A)$ and $\sss(f(A)+\q)$. The
last statement is straightforward. \hfill$\Box$

\bigskip

We already observed  in \cite[Section 5]{dafifoproc} that the kind of results as in
the previous proposition has a moderate interest,  because  the
Krull dimension of $\da$ is compared  to the Krull dimension of $f(A)+\q$, which is  not easy to evaluate
(moreover, if $f^{-1}(\q)=\{0\}$, we have $\da\cong f(A)+\q$ (Proposition \ref{inizio}(3))).

An easy case for  evaluating $\dim(\da)$ is the following.


\begin{prop}\label{precisa}
Let $f: A\rightarrow B$,  $\q$,  and   $\da$ be as in Proposition \ref{inizio}.
Let $f_{\diamond}: A \rightarrow B_{\diamond}:=
f(A) +J$  the ring homomorphism induced from $f$. If we assume that $f_{\diamond}$ is integral (e.g., $f$ is
integral),   then $\dim(\da)=\dim(A)$.
\end{prop}
\noindent
\textsc{Proof}.  By Lemma \ref{integralegamma} and \cite[Theorem 48]{k},
it follows immediately  that $\dim(\da)=\dim(\Gamma(f))=\dim(A)$.\hfill$\Box$

\bigskip

We proceed our investigation  looking for upper and lower bounds of the Krull dimension of $\da$. By Proposition
\ref{spec}, we  know that $\sss(\da)= X \cup U$, where $X:=\sss(A)$ and $U:=\sss(B)\setminus V(J)$ (for the sake
of simplicity, we identify $X$ and $U$ with their homeomorphic images in $\sss(\da)$). Furthermore, again from
Proposition \ref{spec}, we deduce that ideals of the form $\overline Q^{_{_{ f}}}$ can be contained in ideals
of the form $P^{\prime_f}$,  but not vice versa. Therefore,
chains in $\sss(\da)$ are obtained by juxtaposition of two types of chains, one from $U$ ``on the bottom'' and the
other one from $X$ ``on the top'' (where either one or the other may be empty or a single element). It follows
immediately that both $\dim(X)=\dim(A)$ and $\dim(U)$ are lower bounds for $\dim(\da)$ and $\dim (A)+\dim (U)+1$
is an upper bound for $\dim(\da)$  (where, conventionally, we  set $\dim(\emptyset)=-1$).

\begin{oss}\label{attesafrenetica}
Assume that  $J\subseteq {\rm Jac}(B)$. By Proposition \ref{spec}(5), we get that $U$ does not
contain maximal elements of $\sss(\da)$. Hence, in this case, $1+\dim (U)\leq \dim (\da)$.
\end{oss}

Let us define the following subset of $U$:
$$
\mathcal Y_{_{(f,J)}}:=\left\{Q\in U \mid f^{-1}(Q+J)=\{0\}\right\};
$$
it is obvious that $ \mathcal Y_{_{(f,J)}}$ is {\it stable under gene\-rizations},  i.e.,  $Q\in \mathcal
Y_{_{(f,J)}}$, $Q'\in \sss(B)$ and $Q'\subseteq Q$ imply  $Q'\in \mathcal Y_{_{(f,J)}}$. Hence $\dim(\mathcal
Y_{_{(f,J)}})=\ \hskip -4pt \sup\!{\boldsymbol{\{}}\hgt_B(Q)\mid Q\in \mathcal Y_{_{(f,J)}}{\boldsymbol{\}}}$ and
we will denote this integer by $\delta_{_{(f, J)}}$.

\begin{prop} \label{lower} Let $f: A\rightarrow B$,  $\q$,  and   $\da$ be as in  Proposition \ref{inizio};
let $U=\sss(B)\setminus V(J)$ and $\delta_{_{(f, J)}}=
\dim(\mathcal Y_{_{(f,J)}})$ .

\begin{enumerate}
\item[\rm (1)] Let $Q \in \sss(B)$, then $f^{-1}(Q+J)=\{0\}$ if and only if $\overline{Q}^{_f} \ (= (A\times Q) \cap \da)$
is contained in $J_0 \ (= \{0\} \times J)$.
\item[\rm(2)]
for every $Q \in \mathcal Y_{_{(f,J)}}$, the corresponding prime $\overline Q^{_{_{ f}}}$ of $\da$ is contained in
every prime of the form $P^{\prime_f}$.
\item[\rm (3)]
 $\max\{\dim(U),\dim(A)+  1 + \delta_{_{(f,J)}}\} \leq \dim(\da)\,.$
\end{enumerate}
\end{prop}
\noindent

\textsc{Proof}.  (1) Assume that $f^{-1}(Q +J) = \{0\}$. If $(a, f(a) +j)\in \overline{Q}^{_f} $, with $a\in A$
and $j \in J$, then  $f(a) + j \in Q$,  and so $ a \in f^{-1}(Q +J) = \{0\}$, i.e., $a =0$. Therefore, $(a,
f(a) +j) =(0, j) \in J_0$.  Conversely, if  $a \in f^{-1}(Q +J)$, i.e., $f(a) = q +j$ for some $q \in Q$ and $j
\in J$, then $f(a)-j \in Q$,  and so $(a,  f(a)-j) \in  \overline{Q}^{_f} \subseteq J_0$, thus $a =0$.

(2) By Proposition \ref{spec}(1), we have that every ideal of the form $P^{\prime_f}$ contains $J_0$.
The conclusion follows immediately.

(3) By the observation preceding Remark \ref{attesafrenetica}, it is enough to show that $\dim(A)+  1 +
\delta_{_{(f,J)}}\leq \dim(\da)$. If $\mathcal Y_{_{(f,J)}}=\emptyset$ the statement is obvious. Otherwise, let
$Q_0\subset Q_1\subset\z\subset Q_r$
 be a maximal chain  in $\mathcal Y_{_{(f,J)}}$, thus $r=\delta_{_{(f,J)}}$.
Let $P_0\subset P_1\subset \z \subset P_m$ be a chain realizing $\dim(A)$.
By (2) we obtain that $$
\overline Q_0^{_f}\subset\z\subset \overline Q_r^{_f}  \subset   P_0^{\prime_f}\subset \z\subset P_m^{\prime_f},
$$
is a chain in $\sss(\da)$. \hfill$\Box$

\begin{oss}\label{mattinaattesa}
(a) In the situation of Proposition \ref{lower},  note that, if $J $ is contained in the nilradical of $B$, i.e.,
if $V(J) = \sss(B)$, then $ \delta_{_{(f,J)}} = \dim (U)= -1 $. Therefore,
 Proposition \ref{lower}(3) gives $\dim(A) \leq \dim(\da)$. But, in this (trivial) case, we can say more,
precisely  that $\sss(A)$ is homeomorphic to $\sss(\da)$ (Proposition \ref{spec}) and so  $\dim(A) = \dim(\da)$.
As a matter of fact, with the notation of Propositions \ref{pull} and \ref{inizio},  $ \pi_A^\ast: \sss(A)
\rightarrow \sss(\da)$ is a homeomorphism.

(b) Note that, if  $J\not\subseteq \Jac(B)$,   the inequality $ 1+\dim(U) \leq \dim(\da)$ from Remark
\ref{attesafrenetica} can be false, as the following Example \ref{non-d} will show.

(c) Let $f: A\rightarrow B$,  $\q$,  and   $\da$ be as in Proposition \ref{inizio}. If we assume that  $J  \neq
\{0\}$ and that $\da$ and $B$ are integral domains, then, by \cite[Proposition 5.2]{dafifoproc},  $f^{-1}(J) =
\{0\}$ and the subset $\mathcal{Y}_{_{(f,J)}}$ of $\sss(B)$, defined in the previous proposition, is nonempty,
since $(0) \in  \mathcal{Y}_{_{(f,J)}} $, and so $\delta_{_{(f,J)}} \geq 0$.
The following Example \ref{delta}  will show that  $\delta_{_{(f,J)}}$ may be arbitrarily large. Note that
$\delta_{_{(f,J)}}$ may be equal to $-1$ even if $J  \neq \{0\}$, $f^{-1}(J) = \{0\}$, but $B$ is not an integral
~domain. It is sufficient to take $B$ equal to a local zero-dimensional ring not a field, $J$ equal to its maximal
ideal, $A$ any subring of $B$ such that $J \cap A = (0)$, and $f$ be the natural embedding of $A$ in $B$ (e.g.,
$B:= K[X]/(X^2)$, where $K$ is a field and $X$ an indeterminate over $K$, and $A$ any domain contained in $K$). In
this case, $\sss(B) = V(J)$ and so $ \delta_{_{(f,J)}} = -1$.


(d)  Note that, in the situation of Proposition \ref{lower}(1), we can have $\overline{Q}^{_f} \subseteq J_0 \ (=
\{0\} \times J)$ with $Q \supsetneq J$. For instance let $A:=K$,  $B:= K[X, Y]$, $Q := (X, Y)B$, $J:= XB$, and let
$f: A=K \hookrightarrow K[X,Y]=B$ be the natural embedding, where $K$ is a field and $X$ and $Y$ two
indeterminates over $K$.  In this case, $\da \cong  A+ J = K+XK[X, Y]$ (Proposition \ref{inizio}(3)). Clearly,
$f^{-1}(Q) = f^{-1}(Q+J) = f^{-1}(J)= \{0\}$ and  $\overline{Q}^{_f} = J_0 \cong XK[X, Y]$. \end{oss}

\begin{ex} \label{non-d} Let $K$ be a field and $X$ and $Y$ two indeterminates over $K$.
Set $B := K(X)[Y]_{(Y)} \cap K(Y)[X]_{(X)}$.  It is well known that $B$ is a one-dimensional semilocal domain,
having two  maximal ideals $M:= YK(X)[Y]_{(Y)} \cap B$ and $N:= XK(Y)[X]_{(X)} \cap B$. Let $J :=M$, $A:=K$ and
let  $f$ be the natural embedding of $A$ in $B$. Clearly, $f^{-1}(J) = M \cap K = \{0\}$. In this situation, $N
\in \sss(B) \w V(J)$ and so $ \dim (U) =1$. It is easy to see that $\da \cong K +M$ (Proposition \ref{inizio}(3))
is a one-dimensional local domain. Therefore, in this case, we have $2 =1 + \dim(U) > 1 = \dim(\da)$.
\end{ex}

As an immediate consequence of Remark  \ref{attesafrenetica} and Proposition \ref{lower}, we have:

\begin{cor}\label{bestlower}  With the notation of Proposition \ref{lower},
Let $f: A\rightarrow B$,  $\q$,  and   $\da$ be as in Section 2.  If we assume that $J\subseteq {\rm Jac}(B)$ and
that $\delta_{_{(f,J)}}\geq 0$ (e.g., $\da$ and $B$ are integral domains).
Then
$$
1 + \max\!{\boldsymbol{\{}}\dim(A)+  \delta_{_{(f,J)}},\  \dim(U) {\boldsymbol{\}}} \leq \dim(\da)\,.
$$
\vskip -22 pt\hfill $\Box$
\end{cor}

\medskip
The following observations will be useful for Remark
\ref{news2}

\begin{oss}\label{news}
Let $f:A\rightarrow B$, $J$, and $\da$ as in  Proposition \ref{inizio}, and let $Q$ be a prime ideal of $B$.
\begin{enumerate}
\item[(i)] By Proposition \ref{lower}(1), it follows
immediately that $\overline Q^{f}:=(A\times Q)\cap \da\subsetneq J_0:=\{0\}\times J$ if and only if $Q \in
\mathcal  Y_{_{(f,J)}}$ (as defined in Proposition \ref{lower}), i.e., $f^{-1}(Q+J)=\{0\}$ and $Q\nsupseteq J$.
Therefore, $\mathcal Y_{_{(f,J)}}$ is homeomorphic to $\{H \in \sss(\da)\mid H\subsetneq J_0\}$.
\item[(ii)]
If $\da$ and $B$ are integral domains and $J\neq \{0\}$ then, in this situation, $J_0=(0)^{\prime_f}\in\sss(\da)$
and $f^{-1}(J)=\{0\}$  by \cite[Proposition 5.2]{dafifoproc}. Therefore, $Q= (0)\in \mathcal Y_{_{(f,J)}} (\neq
\emptyset)$  and $\overline Q^{f} = f^{-1}(J) \times \{0\} =(0) \subsetneq J_0$; thus, if
$\hgt_{\da}(J_0)<\infty$, $\delta_{_{(f,J)}}(=\dim \mathcal Y_{_{(f,J)}})=\hgt_{\da}(J_0)-1$.

\end{enumerate}
\end{oss}

The next goal is to determine  upper bounds to $\dim(\da)$, possibly sharper than $\dim (A)+\dim (U)+1$.

\begin{thm}\label{carinissima}
 Let $f: A\rightarrow B$,  $\q$,  and   $\da$ be as in  Proposition \ref{inizio}.
 With the notation of Proposition \ref{lower}, assume that $\da$ has finite Krull dimension.  Then
$$
\begin{array}{rl}
&\dim(\da)
 \leq \max \! {\boldsymbol{\{}} \hskip -2pt\dim(A), \   \dim\!\left({A}/{f^{-1}(\q)} \right)  +
 \min\{ \dim(B), \ 1+ \dim (U)  \}  {\boldsymbol{\}}} \\
   & \hskip 10pt \leq \min\!{\boldsymbol{\{}}\!\dim(A)+\dim(U)+1,\ \max\{\dim(A), \ \dim\!\left({A}/{f^{-1}(\q)} \right)  +
  \dim(B) \}  {\boldsymbol{\}}}.

\end{array}
$$
%
\end{thm}

\noindent \textsc{Proof}.   We can  assume that $\sss(B) \neq V(\q)$, because otherwise we already know that
$\dim(\da) = \dim(A)$ (Remark \ref{mattinaattesa}(a)) and so the inequalities hold.

 Let $H_0\subset H_1\subset\z\subset H_n$ be a chain of prime ideals of $\da$ realizing $\dim(\da)$.
 Two extreme cases are possible.

 ~(1) If $H_0\supseteq\{0\}\times\q$ then, by Proposition \ref{spec}(1),
 the chain $H_0\subset H_1\subset \z\subset H_n$ induces  a chain of prime ideals of $A$ of length $n$.
 From  Proposition \ref{lower}(2), we conclude that $\dim(\da)=\dim(A)$.

 ~(2) If $H_n \nsupseteq\{0\}\times\q$. From  Proposition \ref{spec}(2),
 the chain $H_0\subset H_1\subset \z\subset H_n$ induces  a chain of prime ideals of $U$ of length $n$.
 From Proposition \ref{lower}(2), we conclude that $\dim(\da)=  \sup\{\hgt(Q) \mid Q \in U\}= \dim(U)$.

We now consider the general case.

(3) Let $t$ be the maximum index such that $H_t\nsupseteq\{0\}\times\q$, with $0 \leq t \lneq n$. According to the
notations of Proposition \ref{spec}, rewrite the given  chain as follows:
$$
\overline Q_0^{_f}\subset\overline Q_1^{_f}\subset\z\subset\overline Q_t^{_f}\subset\p_{t+1}^{\prime_f}\subset
\p_{t+2}^{\prime_f}\subset\z\subset\p_n^{\prime_f},
$$
where $ Q_0 \subset Q_1 \subset\z\subset  Q_t$ is an increasing chain of prime ideals of $B$, with $Q_t
\not\supseteq J$ (Proposition \ref{spec}(2)), and $ \p_{t+1} \subset \p_{t+2} \subset\z\subset\p_n$ is an
increasing chain of prime ideals of $A$ (Proposition \ref{spec}(1)). Furthermore, by  Proposition \ref{dimp}, we
can find a prime ideal $Q$ in $V(J) \ ( \subseteq \sss(B))$ such that the prime ideal $H_{t+1} =
\p_{t+1}^{\prime_f}$ coincides also with the restriction to $\da$ of the prime ideal $A \times Q$ of $A\times B$,
i.e.,  $H_{t+1} = \p_{t+1}^{\prime_f}= \overline Q^{_f}$. It follows immediately that $P_k\in V(f^{-1}(J))$, for
$t+1\leq k\leq n$.  Therefore,   $  \dim(\da) = (1+t) + (n- t -1)$ with $1+t  \leq \min\{ 1+ \dim(U),\
\dim(B)\}$ and  $n-t -1 \leq  \dim\left({A}/{f^{-1}(\q)} \right)$.

Finally, it is obvious  that $\min\{ \dim(B), \ 1+ \dim(U)  \} \leq \dim(B)$ and that $
\dim\!\left({A}/{f^{-1}(\q)} \right)  +
 \min\{ \dim(B), \ 1+ \dim (U)  \}\leq\dim(A)+\dim(U)+1$.
%
%
\hfill$\Box$


\begin{ex} \label{delta} Let $V$ be a valuation domain with maximal
ideal $\mathfrak M$ such that $V = K + \mathfrak M$, where $K$ is a field isomorphic to the residue field
$V/\mathfrak M$. Let $D$ be an integral domain  with quotient field $K$, and set $B:=D+\mathfrak M$. Assume that
$\dim(V) = n \geq 1$  and that $\mathfrak Q$ is a prime ideal of $V$ with $\hgt_V(\mathfrak Q) = t +1 $, ~$n \geq
t +1 \geq 0$. Set $J := \mathfrak Q \cap B$. By the well known properties of the ``$D+\mathfrak M$
constructions'', $B_{\mathfrak M} = V$ \cite[Exercise 13(1), page 203]{g},  so $J$ is a prime ideal of $B$ and
$\hgt_B(J) = t +1 $.  More precisely, if $(0) \subset {\mathfrak Q}_1 \subset {\mathfrak Q}_2  \subset  \dots
\subset {\mathfrak Q}_t \subset {\mathfrak Q}_{t+1} = \mathfrak Q $  is the chain of prime ideals of $V$ realizing
the height of $\mathfrak Q$, then $Q_0 :=(0) \subset   Q_1:={\mathfrak Q}_1\cap B \subset Q_2 := {\mathfrak Q}_2
\cap B  \subset  \dots  \subset  Q_t:= {\mathfrak Q}_t  \cap B  \subset  Q_{t+1}:= {\mathfrak Q}_{t +1}  \cap B =
J $.  Set $A:=D$ and let $f:A= D \hookrightarrow D+\mathfrak M =B$ be the canonical embedding.  Clearly,
$f^{-1}(J) = \{0\}$ and so it is easy to verify that, in the present situation,
$$
\begin{array}{rl}
\mathcal{Y}_{_{(f,J)}}  :=& \hskip -4pt \{Q\in\sss(B)\mid Q\notin V(J),f^{-1}(Q+J)=\{0\}\}  \\
=& \hskip -4pt  \{Q_k \mid 0 \leq k \leq t \} =  \sss(B)\w V(J)=U
\end{array}
 $$
(see also \cite[Exercise 12(1), page 202]{g}). Therefore,  $ \delta_{_{(f,J)}} =  t = \dim (U)$.   Moreover,
if $ m:= \dim(D) \ (= \dim(A))$  then,  again by the well known properties of the ``$D+\mathfrak M$
constructions'',  $\dim(B) = m +n$ \cite[Exercise 12(4), page 203]{g}.
 Henceforth, in the present example, we have
$ \max\!{\boldsymbol{\{}}\dim(A)+ 1+ \delta_{_{(f,J)}},\   1+ \dim(U)
{\boldsymbol{\}}}  = \dim(A)+ 1+  \delta_{_{(f,J)}} = m+ 1+  t $.

On the other hand,  since $f^{-1}(J) = \{0\}$,  clearly $A/f^{-1}(J)  = A$ and so
$ \max\! {\boldsymbol{\{}} \hskip -2pt\dim(A), \   \dim\left({A}/{f^{-1}(\q)} \right)  +
 \min\{ \dim(B), \ 1+ \dim(U) \} {\boldsymbol{\}}} \ = \ \dim({A}) \  + \\ $
$\min\{ \dim(B), \ 1+ \dim(U)  \}  =  m +\min\{ m+n, \ 1+ t  \}$.
Since $n \geq t+1$,
  then   $\min\{ m+  n, \ 1+ t  \} = 1+t$.
Furthermore, by the fact that $f^{-1}(J) = \{0\}$,  we have $\da \cong A+ J = D+J$ (Proposition \ref{inizio}(3)).
Therefore, from  Corollary \ref{attesafrenetica}  and Theorem \ref{carinissima}, we deduce that $\dim(D+J) = m +
1+ t$.



%
%
%

\end{ex}

\medskip

Let $A\subset B$ be an arbitrary ring extension. We will apply the previous results to the polynomial rings of the
form $A+XB[X]$ and we will show that the bounds given by Fontana, Izelgue and Kabbaj \cite[Theorem 2.1]{fik} in
the very special case where $B$ and $A$ are integral domains coincide to the bounds obtained specializing the general
setting of amalgamated algebras.


%
\begin{oss}\label{significativo}
 Recall that, by \cite[Example 2.5]{dafifoproc}, the ring $A+XB[X]$ (respectively, $A+XB[\![X]\!]$) is
naturally isomorphic to $A\!\Join^{\sigma '}\!\!\!XB[X]$ (respectively, $A\!\Join^{\sigma''}\!\!\!XB[\![X]\!]$),
where $\sigma '$ (respectively, $\sigma ''$) is the inclusion of $A$ into $B':=B[X]$ (respectively, into $B'':= B[\![X]\!]$).
\end{oss}
\begin{cor}\label{poli1}
Let $A\subseteq B$ be a ring extension  and $X$ an indeterminate over $B$. Set
$$
 \delta'_{(A,B)}:=\sup\!{\boldsymbol{\{}}\hgt_{B[X]}(Q)\mid Q\in\sss(B[X]), X\notin Q,(Q+XB[X])\cap A=\{0\}{\boldsymbol{\}}}.
$$
Then
$$
\begin{array}{rl}

\max\!{\boldsymbol{\{}}\!\dim(A) +  1+ \delta'_{(A,B)},\  & {\hskip -9pt}\dim(B[X,X^{-1}]){\boldsymbol{\}}}
 \leq \dim(A+XB[X])\leq \\
\leq  & {\hskip -6pt}
\dim(A)+\dim(B[X])\,.
 \end{array}
 $$
\end{cor}

\noindent \textsc{Proof}. Let $B' := B[X]$ and $J' := XB[X]$.  As observed above (Remark \ref {significativo}), we
know that $A\!\Join^{\sigma '}\!\!\!J' =A+XB[X]$. From the definitions, it is easy to see that
$\delta_{_{(\sigma', J')}}=\delta'_{(A,B)}$. Moreover, since $\dim(B[X,X^{-1}])=\sup\{\hgt_{B[X]}(Q)\mid
Q\in\sss(B[X]), X\notin Q\} = \dim(U)$ (where $U$, in this case, is homeomorphic to $\sss(B[X])\setminus V(J')$)
and $\sigma'^{-1}(J')=A\cap XB[X]=\{0\}$, the conclusion follows from Proposition \ref{lower}(3) and Theorem
\ref{carinissima}.\hfill$\Box$

%
%
%

\begin{oss}\label{news2}
Let $A\subseteq B$  integral domains and and let  $N:=A\w\{0\}$. In \cite[Theorem 2.1]{fik}, Fontana, Izelgue
and Kabbaj  proved that
$$
\begin{array}{rl}
\max\!{\boldsymbol{\{}}\!\dim(A) + \dim(N^{-1}B[X]), & {\hskip -8pt} \dim(B[X]){\boldsymbol{\}}}
\leq \dim(A+XB[X])\leq \\
\leq  & {\hskip -6pt}
\dim(A)+\dim(B[X])\,.
 \end{array}
 $$
By  \cite[Theorem 1.2(a) and Lemma 1.3]{fik},  we know that
$$
\dim(N^{-1}B[X]) = \hgt_{A+XB[X]}(XB[X]) = 1 + \lambda'_{(A,B)},
$$
 where
$$
\lambda'_{(A,B)}:=\sup\!{\boldsymbol{\{}}\!\dim\left(B[X]_{{\boldsymbol q}[X]}\right) \mid {\boldsymbol q}\in\sss(B), \   {\boldsymbol q} \cap A= (0){\boldsymbol{\}}}\,.
$$
 From Remark \ref{news}(iii)
 and the proof of Corollary \ref{poli1},
 we deduce the equality
$\hgt_{A+XB[X]}(XB[X])=1+\delta'_{(A,B)}=1+\lambda'_{(A,B)}$, hence  $\delta'_{(A,B)}=\lambda'_{(A,B)}$; moreover,
we have  $\dim B[X]=\dim B[X,X^{-1}]$, by  \cite[Proposition 1.14]{abdfk}. Therefore, in particular, we reobtain
Fontana, Izelgue and Kabbaj's result on the dimension of the integral domain $A+XB[X]$. This fact provides further
evidence on the sharpness of the bounds obtained in Proposition \ref{lower}(3) and Theorem \ref{carinissima}, in
the general setting of amalgamated algebras.

 %
\end{oss}
\medskip

We consider now the case of power series rings of the type $A+XB[\![X]\!]$ for  arbitrary ring extensions  $A\subset B$.
%
\begin{cor}\label{sipuodiremeglio}
Let $A\subset B$ be a ring extension and $X$ an indeterminate over $B$. Set
$$
\delta^{''}_{(A,B)}:=\sup\!{\boldsymbol{\{ }}\hgt_{B[\![X]\!]}(Q)\mid  Q\in\sss(B[\![X]\!])\w V(X),
(Q+XB[\![X]\!])\cap A=\{0\}{\boldsymbol{ \} }}.
$$

Then
$$
\begin{array}{rl}

\max\!{\boldsymbol{\{}} \dim(A) +  1+ \delta^{''}_{(A,B)},\ 1 + & {\hskip -9pt}\dim(B[\![X]\!][X^{-1}]){\boldsymbol{\}}}
\leq\dim(A+XB[\![X]\!])\leq \\
\leq  & {\hskip -6pt}
 1+ \dim(A) + \dim(B[\![X]\!][X^{-1}]).
 \end{array}
$$
\end{cor}

\noindent
\textsc{Proof}. Keeping in mind the statements and the notation of  Remark \ref{significativo}, it follows immediately that  $\delta_{(\sigma'', XB[\![X]\!])}=\delta^{''}_{(A,B)}$. Moreover,
recalling that $U$, in this case, is homeomorphic to $\sss(B[\![X]\!])\setminus V(X)$,
it
is easy to see that $\dim(U)=\dim(B[\![X]\!][X^{-1}])$.
Finally, note that $\min\{\dim(B[\![X]\!]), \  1+\dim (U)\} = 1+\dim (U)$, since every
maximal ideal of $B[\![X]\!]$ contains $X$ \cite[Chapter 1, Exercise 5(iv)]{am}.
 The conclusion is now a straightforward consequence of   Corollary \ref{attesafrenetica}
  and Theorem \ref{carinissima}.
\hfill$\Box$

%

\begin{oss}

By applying Corollary \ref{sipuodiremeglio} and Remark \ref{mattinaattesa}, it follows that, if $B$ is an integral
domain, then
$$
\begin{array}{rl}
1+\max\!{\boldsymbol{\{}}\dim(A) +\delta^{''}_{(A,B)},\ & \hskip -9 pt
\dim(B[\![X]\!][X^{-1}]){\boldsymbol{\}}}\leq\dim(A+XB[\![X]\!])\leq \\
\leq & \hskip -6 pt
 1+\dim(A)+\dim(B[\![X]\!][X^{-1}]).
 \end{array}
$$

Now, we can compare our lower bound with that given by Dobbs and Khalis's Theorem (\cite[Theorem 11]{doka}).
Setting
$$\lambda^{''}_{(A,B)}:=\sup\!{\boldsymbol{\{}}
\dim\!\left(B[\![X]\!]_{{\boldsymbol q}[\![X]\!]}\right) \mid {\boldsymbol q}\in\sss(B),
 \ {\boldsymbol q}\cap A=(0){\boldsymbol{\}}}\,,
$$
they prove that
$$
\begin{array}{rl}
1+\max\!{\boldsymbol{\{}}\dim(A) +\lambda^{''}_{(A,B)},\ & \hskip -9 pt
\dim(B[\![X]\!][X^{-1}]){\boldsymbol{\}}}\leq\dim(A+XB[\![X]\!])\leq \\
\leq & \hskip -6 pt
 1+\dim(A)+\dim(B[\![X]\!][X^{-1}]).
 \end{array}
 $$
It is clear that $ \dim\!\left(B[\![X]\!]_{{\boldsymbol q}[\![X]\!]}\right)= \hgt_{B[\![X]\!]}({\boldsymbol
q}[\![X]\!])$. Moreover, it is immediately seen that, if ${\boldsymbol q}\in\sss(B)$ and ${\boldsymbol q}\cap A=
(0)$, then $({\boldsymbol q}[\![X]\!]+XB[\![X]\!])\cap A= (0)$.   Since the set $ {\boldsymbol{\{}} {\boldsymbol
q}[\![X]\!] \in \sss(B[\![X]\!]) \mid {\boldsymbol q} \in \sss(B) \mbox{ and } {\boldsymbol q}\cap A=\{0\}
{\boldsymbol{\}}}$ is a subset of ${\boldsymbol{\{}}Q\in \sss(B[\![X]\!]) \mid X \notin Q \mbox{ and }  (Q
+XB[\![X]\!])\cap A=\{0\}{\boldsymbol{\}}} $, we have $\lambda^{''}_{(A,B)}\leq \delta^{''}_{(A,B)}$. It is
natural to ask, as in the polynomial case: does $\lambda^{''}_{(A,B)} =\delta^{''}_{(A,B)}$ hold? At the moment,
the answer to this question is open. However, by \cite[Theorem 7]{doka},  we observe that the answer could be
negative if
$$ \hgt_{A+XB[\![X]\!]}(XB[\![X]\!]) = 1 + \delta^{''}_{(A,B)}$$
and  $\lambda^{''}_{(A,B)} \lneq
   \sup\{\hgt_{B[\![X]\!]}(Q) \mid Q \in{\bf{\Lambda}}_{(A,B)} \}\,,$
where ${\bf{\Lambda}}_{(A,B)}$, as in \cite[Theorem 7]{doka}, is defined to be ${\bf{\Lambda}}_{(A,B)}=\{ Q \in
\sss(B[\![X]\!])\mid X\notin Q, \ Q \subset (\boldsymbol q, X),\ \mbox{for some}\ \boldsymbol q $ $\in \sss(B) \mbox{
with } \boldsymbol q \cap A = (0) \}$.
\end{oss}

\begin{ex}
It is possible to construct an infinite dimensional ring of the type $\da$, where $A$ is a finite dimensional
ring. In this situation, $B$ must be  a infinite dimensional ring by Theorem \ref{carinissima}. For instance, let
$A:=\mathbb C$ be the field of complex numbers, let $Y$ be an indeterminate over $\mathbb C$, and let $R:=\mathbb
C[\{Y^{1/n}\mid n\in\nat\w\{0\}\}]$. Consider the maximal ideal $\mathfrak M$ of $R$ generated by the set
$\{Y^{1/n}\mid n\in\nat\w\{0\}\}$. Set $B:=R_{\mathfrak M}$, and consider the ring $A+XB[\![X]\!] \ (\cong
A\Join^{\sigma''} XB[\![X]\!]$, according to notation of Remark \ref{significativo}). Then, by \cite[Example
3]{doka},  $B$ is a one-dimensional non-discrete valuation domain and  $\hgt_{A+XB[\![X]\!]}(XB[\![X]\!])=\infty$,
and thus $\dim(A+XB[\![X]\!])=\infty$.
\end{ex}
\smallskip

The next two examples show that the upper bound and lower bound of Theorem \ref{carinissima} and Proposition
\ref{lower}(3)  are ``sharp", in the sense that $\dim(\da)$ may be equal to each of the  two  numerical terms
appearing in the  first  inequality
(respectively,  in the inequality)  of Theorem \ref{carinissima} (respectively, Proposition \ref{lower}(3)).

\begin{ex} Let $A$ be a valuation domain such that $\dim(A)=n \geq 3$,
let $\{0\}\subset \p_1\subset\p_2\subset \dots \subset \p_n$
be a chain of prime ideals of $A$ realizing $\dim(A)$, and let
$x_h \in \p_{h+1}\w\p_h$, with $1 \leq h \leq n-2$ and $(x_h) \neq P_{h+1}$.
Since $A$ is a valuation domain, it is easily seen that $V(x_h)=V(\p_{h+1})$, and thus
$\dim(A/(x_h))=\dim(A/\p_{h+1})=n -(h+1)$. Set $B:= A/(x_h)$, $f:A \twoheadrightarrow B$ the canonical projection,
$Q_k := P_{k}/(x_h)$ for $h+1 \leq k \leq n$, and $J := Q_{h+j}$  for some $1 \leq j \leq n-h$. In this case, by
Proposition \ref{su}, $\dim(\da) = \dim(A \times B) = \max\{\dim(A), \dim(B)\} = \dim(A)=n$. Note also that $
\dim\left({A}/{f^{-1}(\q)}\right) = n -(h+j) \leq  n -(h+1) =\dim(B)$, and
$$
 \dim (U) =
\left\{
\begin{array}{l}
-1 \mbox{\,, \; if $ j =1$,} \\
 j -2   \mbox{\,, \; if $ 1 <  j \leq n -h$} .
\end{array}
\right.
$$
 It is also easy to see that $f^{-1}(Q+J) \neq \{0\}$
for all $Q \in \sss(B)$ and so in this case $\delta_{_{(f,J)}} =-1$, for all $1 \leq j \leq n-h$.
 Moreover, in the present situation, $\da$
is a local ring, but it is not an integral domain since $f^{-1}(J) \neq \{0\}$
(see \cite[Proposition 5.2]{dafifoproc}).

Consider now a chain
$H_0\subset H_1\subset\z\subset H_n$
of prime ideals of $\da$ realizing $\dim(\da)$.  Two cases are possible.

~$\bullet$  If $1 \lneq j \leq n-h$, then the previous chain  (realizing $\dim(\da)$) is of the type:
$$
\begin{array}{rl}
((0)\neq) P_0^{\prime_f}\subset & \hskip -5pt P_1^{\prime_f}  \subset\z
\subset   P_h^{\prime_f} \subset \\
\subset & \hskip -5pt  P_{h+1}^{\prime_f} =\overline{Q}_{h+1}^{_f}  \subset\z
\subset P_{h+j -1 }^{\prime_f} =\overline{Q}_{h+j -1}^{_f} \subset \\
& \hskip -5pt \subset P_{h+j}^{\prime_f}  \subset\z \subset P_n^{\prime_f}
\end{array}
$$
 (where  $P_{k}^{\prime_f} = \overline{Q}_{k}^{_f}$ also for $h+j \leq k \leq n$, but in this case $Q_k \supseteq J$);

$\bullet$   If   $ j=1$, then the previous chain realizing $\dim(\da)$ is of the type:
 $$
 ((0)\neq) P_0^{\prime_f}\subset P_1^{\prime_f}\subset\z\subset P_h^{\prime_f} \subset\z\subset  P_n^{\prime_f}$$
 and none of the $P_{k}^{\prime_f} $ is equal to a  $\overline{Q}_{k}^{_f} $ for $Q_k \not\supseteq J$.

In the present example, the inequality of Corollary \ref{attesafrenetica} gives  back  the inequality $
\max\!{\boldsymbol{\{}}\dim(A)+ 1+ \delta_{_{(f,J)}},\ 1+ \dim(U) {\boldsymbol{\}}} $ = $
\max\!{\boldsymbol{\{}}n+ 1 -1,\ 1+  (j-2) {\boldsymbol{\}}}$ $\leq$ $n =  \dim(\da)$. The  first inequality of
Theorem \ref{carinissima} gives $\dim(\da) =n   \leq \max \! {\boldsymbol{\{}}\! n, \   n-(h+j) + \min\{n-(h+1), \
1+ (j-2)  \}  {\boldsymbol{\}}}$ \ = \ $ \max \! {\boldsymbol{\{}}\! \dim(A), \;$ $ \dim\left({A}/{f^{-1}(\q)}
\right)  + \min\{ \dim(B), \ 1+ \dim (U)  \}  {\boldsymbol{\}}}$.


\end{ex}

\begin{ex} Let $K$ be a field and let $V$ and $W$ be
two incomparable finite dimensional valuation domains having same field of quotients $F$. Assume that $V$ and $W$
are $K$-algebras, that $V = K +\mathfrak M$ and $W = K +\mathfrak N$ where $\mathfrak M$ (respectively, $\mathfrak
N$) is the maximal ideal of $V$ (respectively, $W$), and that $\dim(V)= m\geq 1$ and $\dim(W) = n\geq 1 $.  Set $T
:= V \cap W $. It is well known that  $T$ is a finite dimensional B\'ezout domain with quotient field $F$ and with
two maximal ideals $M: = \mathfrak M \cap T$ and $N:= \mathfrak N \cap T$ such that  $T_M= V$ and $T_N = W$, and
so $\dim(T) =\max\{m, n\}$  \cite[Theorem 101]{k}.
 Let $D$ be an integral domain of Krull dimension $d$ with quotient field $K$.
 Since $D$ is embedded naturally in $V \ (= K+\mathfrak M)$ and $W \ (= K +\mathfrak N)$,
we have also a natural embedding $\iota: D \hookrightarrow T $.

In this situation, using the standard notation of
the $\da$ construction, when $A := D$, $B:= T$, $J := M$, and $f:= \iota$, we have that   the ring $D + M$
(subring of $T$) is canonically isomorphic to $D\!\Join^{\iota}\!\! M$, by \cite[Example 2.6]{dafifoproc}.
Moreover, $f^{-1}(J) = M \cap D =\{0\}$ and so
   $\dim(A/f^{-1}(J))=\dim(D) = d$, and
$\dim(U) = \max\{ m-1, n\}\,.$

 It is easy to verify that if $(0)= \mathfrak Q_0 \subset \mathfrak Q_1 \subset ...\subset \mathfrak Q_m = \mathfrak M$
 are the prime ideals of $V$, then $\{ Q_k := \mathfrak Q_k \cap B  \mid 0\leq k \leq m-1\}$  coincides with the ~set
  $ \{ Q \in \sss(B) \w V(J) \mid f^{-1}(Q+J) = (Q+J) \cap D =\{0\} \}$.
  Therefore, $\delta_{_{(f,J)}} =m-1$.
On the other hand, it is easy to verify that
$
\dim(D+M) = \max \{ m+d, n \}$.

 In the present example,  the inequality of Proposition \ref{lower}(3) gives  back  the inequality
$ \max\!{\boldsymbol{\{}}\dim(A)+ 1+ \delta_{_{(f,J)}} ,\  \dim(U) {\boldsymbol{\}}} $ = $
\max\!{\boldsymbol{\{}}d + 1 + m-1  ,\ \max\{m-1, n\}{\boldsymbol{\}}}$ $\leq$ $\max\{m+d, n\} = \dim(\da)=
\dim(D+M)$.  Therefore, if $n> m+d$, then $\dim(\da) = \dim(U)$. By the first inequality of Theorem
\ref{carinissima}, it follows that $\dim(\da) = \max\{m+d, n\}  \leq \max \! {\boldsymbol{\{}}\! d, \  d +
\min\{\max\{m,n\}, \ 1+ \max\{m-1,n\}  \}  {\boldsymbol{\}}}$ = $ \max \! {\boldsymbol{\{}}\! \dim(A),$ $
\dim\left({A}/{f^{-1}(\q)} \right)  + \min\{ \dim(B), \ 1+\dim(U)  \}  {\boldsymbol{\}}}$. Therefore, if $m +d
\leq n$, then $ n = \dim(D+M)  = \dim(D\!\Join^{\iota}\!\! M) = d + \min\!{\boldsymbol{\{}}\!\max\{m,n\}, \ 1+
\max\{m-1,n\}{\boldsymbol{\}}}$.

\end{ex}

\noindent \textsc{Acknowledgements}. The authors are very grateful to the referee for several suggestions and
comments that greatly improved the paper.


\end{document}